# Weak systems of determinacy and arithmetical quasi-inductive definitions.

P.D.Welch[1]

School of Mathematics
University of Bristol
England

*May 18, 2009*

### Abstract

We locate winning strategies for various $\Sigma_3^0$-games in the $L$-hierarchy in order to prove the following:

**Theorem 1.** KP + $\Sigma_2$-Separation $\vdash$ $\exists \alpha$ $L_\alpha \vDash$ "$\Sigma_2$-KP + $\boldsymbol{\Sigma_3^0}$-Determinacy."

Alternatively: $\boldsymbol{\Pi_3^1\text{-}CA_0} \vdash$ "there is a $\beta$-model of $\boldsymbol{\Delta_3^1\text{-}CA_0} + \boldsymbol{\Sigma_3^0}$-Determinacy." The implication is not reversible. (In fact the antecedent here may be replaced with $\Pi_2^1(\Pi_3^1)$-CA$_0$: $\Pi_2^1$ instances of Comprehension with only $\Pi_3^1$-lightface definable parameters.)

**Theorem 2.** KP + $\Delta_2$-Comprehension + $\Sigma_2$-Collection + **AQI** $\nvdash \Sigma_3^0$-Determinacy.

(Here **AQI** is the assertion that every arithmetical quasi-inductive definition converges.) Alternatively:

$\boldsymbol{\Delta_3^1\text{-}CA_0}$ + **AQI** $\nvdash \Sigma_3^0$-Determinacy.

Hence the theories: $\boldsymbol{\Pi_3^1\text{-}CA_0}$, $\boldsymbol{\Delta_3^1\text{-}CA_0} + \Sigma_3^0$-Det, $\boldsymbol{\Delta_3^1\text{-}CA_0} + \mathbf{AQI}$, and $\boldsymbol{\Delta_3^1\text{-}CA_0}$ are in strictly descending order of strength.

# 1 Introduction

The work in this paper was, initially at least, motivated by trying to see how the theory of *arithmetical quasi-inductive definitions* (**AQI** as defined below) fits in with other subsystems of second order number theory. We had been working with one example of such a definition, essentially a *recursive* quasi-inductive definition, and had calculated certain ordinals where such definitions reached fixed points or exhibited a looping convergence [17] and [16]. Earlier J. Burgess [2] had in fact distilled from H. Herzberger's notion of a *revision sequence* [8] the notion of arithmetically quasi-inductive, and shown that the same ordinals appeared. (Herzberger's notion was connected with a ''truth operator'' and thus, strictly speaking is not arithmetic, but just beyond; however this only makes for a trivial difference.) Other examples of constructions involving such quasi-inductive definitions have appeared in the theory of truth [4], and in theoretical computer science: S.Kreutzer in [10] uses essentially arithmetical quasi-inductive definitions to formulate a notion of semantics for partial fixed point logics over structures with infinite domains in order to separate away this logic from inflationary fixed point logic.

Here however we rather mention some of the possibilities that connect these concepts with potential proof theoretical results on the way to looking at ordinal systems for $\Pi_3^1$-CA$_0$: for $\Pi_2^1$-CA$_0$, by work of Rathjen [13],[12], we have that this second level of Comprehension is tied up with the theory of arbitrarily long finite $\Sigma_1$-elementary chains through the $L_\alpha$-hierarchy: indeed the first level $L_\alpha$ which is an infinite tower of such models, is the first whose reals form a $\beta$-model of $\Pi_2^1$-CA$_0$. The same occurs for $\Pi_3^1$-CA$_0$: the first $L_\gamma$ whose reals form a $\beta$-model of $\Pi_3^1$-CA$_0$ is the union of an infinite tower of models $L_{\zeta_n} \prec_{\Sigma_2} L_{\zeta_{n+1}}$ and presumably one will need to analyse finite chains of such models to get at an ordinal system for this theory. Seeing that **AQI** is connected with levels $L_\zeta$ of the Gödel's $L$-hierarchy, with $\Sigma_2$-end extensions, (albeit only chains of length 1) analysing the proof theoretic strength of **AQI** would be a natural stepping stone.

Other notions of inductive definition have been tied to *determinacy*. Positive monotone arithmetical operators have fixed points bounded by the first admissible ordinal $\omega_1^{ck}$, and in turn strategies for recursive open (that is $\Sigma_1^0$) games are either in $L_{\omega_1^{ck}}$ (for player $I$) or definable over it (for the ''closed'' player $II$). Solovay [12] showed that, remarkably, for $\Sigma_2^0$-games, strategies for player $I$ in such games occur in $L_\sigma$ where $\sigma$ is the closure ordinal of $\Sigma_1^1$ monotone inductive definitions (and for Player $II$ they lie in the next admissible set beyond it). Tanaka [15] formulated a subsytem of analysis related to $\Sigma_1^1$-monotone inductive definitions, $\Sigma_1^1$-$\mathrm{MI}_0$, and showed that over $\mathrm{RCA}_0$ it was equivalent to $\Sigma_2^0$-Determinacy. Our original intention had been to tie in **AQI** with some level of determinacy. Is there anything at all analogous for **AQI**?

Turning naturally to $\Sigma_3^0$-Determinacy, the location of strategies for such games in the constructible $L$-hierarchy seems to be unknown. There is little published on this question  John states some results on this in in [9]. However this does not really yet reveal (at least to us) where such strategies lie. A closer reading of Davis' proof of $\Sigma_3^0$-Determinacy showed that it was provable in $\mathrm{KP} + \Sigma_2$-Separation, and thus that winning strategies appear in the least model which is an infinite tower of the form of a union of a chain of submodels $L_{\zeta_n} \prec_{\Sigma_2} L_{\zeta_{n+1}}$. Whilst we always thought it would be a happy coincidence if $\Sigma_3^0$-Determinacy matched up exactly with **AQI** we never really believed it would be so, and the theorems here show this. However they are *extremely* close. We have not located the exact ordinal where the winning strategies (for either player) appear, but we have shown that for some games they must appear *after* the first $\Sigma_2$-extendible $\xi_0$ - that is $\xi_0$ initiates a $\Sigma_2$-chain of length 1: $L_{\xi_0} \prec_{\Sigma_2} L_{\xi_1}$ (indeed after the first $\Sigma_2$-admissible $\mu$ so that the reals of $L_\mu$ are closed under boldface **AQI**). However for all $\Sigma_3^0$ games they must appear before the first $\gamma$ (indeed are strictly bounded below) where the reals of $L_\gamma$ are a model of ''$\Pi_2^1(\Pi_3^1)$-$\mathrm{CA}_0$'' - that is they are closed under instances of $\Pi_2^1$ Comprehension, with only $\Pi_3^1$ (lightface) definable parameters allowed (a more precise definition of a better bound is below). Such an ordinal $\gamma$ occurs before the beginning of the least $\Sigma_2$-chain of length 2: $L_{\zeta_0} \prec_{\Sigma_2} L_{\zeta_1} \prec_{\Sigma_2} L_{\zeta_2}$. It seemed to us that even with extra effort, finding exactly at which level each $\Sigma_3^0$ game had a strategy might not be very much more illuminating: any such precise characterisation of these levels might well be in terms of $\Sigma_3^0$ games anyway, rather than some exterior defined subsystem of second order number theory, or of particular properties of certain levels of the $L$-hierarchy.

## 1.1 Preliminaries

By $\Sigma_2$-Separation we mean the usual Axiom of Separation but restricted to formulae that are $\Sigma_2$ in the Levy hierarchy. By ''KP'' we shall mean the usual axioms of Kripke-Platek set theory, but we shall assume these include the Axiom of Infinity. By ''$KPI_0$'' we shall mean the conjunction of the axiom of extensionality together with the assertion ''For every set $x$ there is a set $y$ with $x \in y \wedge (KP)_y$ ''. By ''KPI'' we shall mean ''$KP + KPI_0$.'' By ''$\Sigma_2$-KP'' we shall mean KP with the Comprehension and Collection Axiom schemes reinforced to allow for $\Delta_2$ and $\Sigma_2$ formulae respectively.

We note that KP alone does not prove $\Sigma_1$-Separation. The Separation axioms are themselves essentially ''boldface'' axioms, as they allow parameters into the axiom schemes. The class of all admissible ordinals other than $\omega$, we shall denote by ''ADM''. For information on admissibility theory the reader may consult [1]. The reals of a model of $KP + \Sigma_i$-Separation form a model of $\mathbf{\Pi}^1_{i+1}$-$CA_0$ for $i \in \{1, 2\}$. We follow the definitions and development of the theories $\mathbf{\Pi}^1_{i+1}$-$CA_0$ (the latter we take to include the set induction axiom) of [14]. The set of reals belonging to $L_\alpha$, where $\alpha$ is least such that $L_\alpha$ is a model of $KP + \Sigma_i$-Separation, are those of the minimum $\beta$-model of $\mathbf{\Pi}^1_{i+1}$-$CA_0$, $(i > 0)$ (cf. [14],VII.5.17). It is well known, and easy to see, that if $L_\alpha$ is a countable model of $KP + \Sigma_i$-Separation, then $L_\alpha$ is a union of an infinite $\Sigma_i$-elementary chain of submodels $L_{\zeta_k} \prec_{\Sigma_i} L_{\zeta_j}, \cdots \prec_{\Sigma_i} L_\alpha$. The least such $\alpha$ then being the first that is the union of an $\omega$-length such chain. As a consequence $KP + \Sigma_i$-Separation proves the existence of $\beta$-models whose reals code $\Sigma_i$-elementary chains of length, say 2: $L_{\zeta_0} \prec_{\Sigma_i} L_{\zeta_1} \prec_{\Sigma_i} L_{\zeta_2}$.

Our set theoretical notation is standard. If $z = \langle u, v \rangle$ is an ordered pair, then we use the functions $(z)_0 = u$, $(z)_1 = v$ to denote the unpairing functions. For $D$ a class of ordinals we let $D^*$ be the closed class of its limit points.

### 1.1.1 On the constructible hierarchy

We shall use some further facts about the $L$-hierarchy. In particular that if $\omega\alpha = \alpha$ then $L_\alpha$ is equal to the level $J_\alpha$ of the Jensen hierarchy. (We shall write $L_\alpha$ rather than $J_\alpha$ when we know that they are equal, and for the ordinal heights of the models of our theories, this will always be the case.)

We shall write $\Sigma_n(L_\alpha)$ for the relations $\Sigma_n$ (for $n \leq \omega$) over $\langle L_\alpha, \in \rangle$ possibly with paramters from $L_\alpha$. If we wish to display the parameters $p, q, r\ldots$ we shall write $\Sigma_n^{L_\alpha}(\{p, q, r\})$ (and similarly for $J_\alpha$ if needed in place of $L_\alpha$.) We shall not make much use of the fine-structure theory but we shall use the $\Sigma_1$-skolem function $h^1_\alpha(v_0, v_1)$ which is uniformly definable over such levels. Thus if $\langle \varphi_k \mid k < \omega \rangle$ is a recursive listing of the $\Sigma_1$ formulae of $\mathcal{L}_{\dot\in}$, with say $\varphi_k$ with free variables amongst those displayed as $\varphi_k(v_0, \ldots, v_{n(k)})$, and if $L_\alpha \vDash \exists v_0 \varphi_k(v_0, x_1 \ldots, x_{n(k)})$ then $L_\alpha \vDash \varphi_k(h^1_\alpha(k, \langle p^1_\alpha, \vec{x}_i \rangle), x_1 \ldots, x_{n(k)})$ where $p^1_\alpha$ is the (possibly empty) $\Sigma_1$-*standard parameter*. The reader should remark that we are working very low down in the $L$-hierarchy: any level will have at most some $\Sigma_2$ definable map from $\omega$ onto the whole structure; every $L_\alpha$ of interest will have only $\omega$ as the single infinite cardinal, and any stnadard parameters that occur can be taken to be single ordinals. Also for any $\alpha$ there will be always a $\Delta_1(L_\alpha)$ definable map of $\alpha$ onto $L_\alpha$ (although for $\alpha$ not closed under Gödel pairing this is not necessarily a parameter free definition; for those $\alpha$ which are Gödel closed, the definition is uniform in $\alpha$.)

Indeed we are so far down in the $L$-hierarchy that we always have *uniformly $\Sigma_n$-definable skolem functions*. For this recall that $\beta_0$ is defined as the least $\beta$ such that $L_{\beta_0} \vDash \text{ZF}^-$. (It is well known that the reals of $L_{\beta_0}$ are those of the minimal model of $\Pi^1_\infty\text{-CA}_0$, the minimal model of full second order analysis.) Friedman showed in [5] that if $J_\alpha \vDash$ ''$V = \text{HC}$'' (that is every set is hereditarily countable), then for any $n$ there is $g^n_\alpha$ a $\Sigma_n$-Skolem function for $J_\alpha$ $\Sigma_n$-definable without parameters over $J_\alpha$. (The point here is the phrase ''without parameters''.)

**Remark 3.** It is well known that for $n > 1$ there is no uniform $\Sigma_n$-definable $\Sigma_n$-skolem function for all the levels of the $J_\alpha$ hierarchy. It is perhaps less well known that down low, there are in fact such as we now show.

**Theorem 4.** (Uniform $\Sigma_n$-Skolem Functions) *For every $n < \omega$ there is a single $\Sigma_n$-definition of a partial function $h^n$, which defines a $\Sigma_n$-skolem function over any $\langle J_\alpha, \in \rangle$ with $\alpha < \beta_0$.*

**Proof**. It is proven in [5] that for a fixed $\alpha_0$ by using an induction on $n < \omega$ that $\Sigma_n$-skolem functions exists. For the theorem under discussion here, any $J_\alpha$ for $\alpha \le \beta_0$ satisfies the requirement $J_\alpha \vDash$ ''$V = \text{HC}$''. For an $\alpha_0$ as given and starting with $n = 2$ Friedman's definition for a skolem function $g^2_{\alpha_0}$, is in terms of the $\Sigma_1$-skolem function which is indeed uniformly parameter free $\Sigma_1$ over all $\alpha$. An inspection of the definition (on p3328 op.cit.) shows however that it does nor depend on $\alpha_0$ in any way. We thus have that for any $\alpha \le \beta_0$ there is indeed a $\Sigma_2$-skolem function uniformly parameter free $\Sigma_2$-definable over $J_\alpha$, which we shall call $h^2_\alpha$. Turning to $n = 3$ his definition for the fixed $\alpha_0$ is in terms of this $\Sigma_2$-definable skolem function $g = g^2_{\alpha_0}$ (as defined at the bottom of p3328) for $J_{\alpha_0}$, but is otherwise independent of $\alpha_0$. And we have just argued that the $\Sigma_2$-skolem functions $h^2_\alpha$ may be assumed to be here uniformly definable. Hence so will be the $\Sigma_3$-skolem functions, defined using them, that we may call $h^3_\alpha$ for all $\alpha \le \beta_0$. Similarly for other $n < \omega$. Q.E.D.

Consequently, every $x \in J_\alpha$ (for $\alpha < \beta_0$) is of the form $h_\alpha^n(i, n)$ some $i$, $n < \omega$. In fact the argument in the above proof readily works *uniformly* beyond $\beta_0$ up until the least $\alpha + 1$ such that $J_{\alpha+1} \vDash$ ''*there is an uncountable cardinal*.''

We shall also use the fact that if $T$ is a recursive theory of $\Sigma_1$ sentences, and $\alpha$ is the least such that $L_\alpha \vDash \mathrm{KP} + T + V = \mathrm{HC}$ (that every set is hereditarily countable) then every $x \in L_\alpha$ is definable by some $\Sigma_1^{L_\alpha}$ term. This is because there will be a $\Sigma_1^{L_\alpha}$ map, perhaps partial, of $\omega$ onto $L_\alpha$: hence for every $x$ there is $i, n < \omega$ with $x = h_\alpha^1(i, n)$.

## 1.2 Arithmetical quasi-inductive definitions.

Let $\Gamma \colon \mathcal{P}(\omega) \to \mathcal{P}(\omega)$ be any arithmetic operator (that is ''$n \in \Gamma(X)$'' is arithmetic; we emphasise that $\Gamma$ need be neither monotone nor progressive). We define the following iterates of $\Gamma$: $\Gamma_0(X) = X$; $\Gamma_{\alpha+1}(X) = \Gamma(\Gamma_\alpha(X))$; $\Gamma_\lambda(X) = \liminf_{\alpha \to \lambda} \Gamma_\alpha(X) = \cup_{\alpha < \lambda} \cap_{\lambda > \beta > \alpha} \Gamma_\beta(X)$. Following Burgess we say that $Y \subseteq \omega$ is *arithmetically quasi-inductive* if for some such $\Gamma$, $Y$ is (1-1) reducible to $\Gamma_{\mathrm{On}}(\varnothing)$. Any such definition has a least countable $\xi = \xi(\Gamma)$ with $\Gamma_\xi(\varnothing) = \Gamma_{\mathrm{On}}(\varnothing)$. If we let $\zeta$ denote the supremum of all such $\xi(\Gamma)$, then we have:

**Proposition 5.** *(Burgess [2] Sect.14) (i) $\zeta$ is the least $\Sigma_2$-extendible ordinal; that is the least $\zeta$ so that there is a $\Sigma > \zeta$ with $L_\zeta \prec_{\Sigma_2} L_\Sigma$.*
   *(ii) A set $Y$ is arithmetical quasi-inductive iff $Y \in \Sigma_2(L_\zeta)$.*

In general we shall stay with the notation that $\Sigma$ denotes the ordinal height of the least proper $\Sigma_2$-end extension of $L_\zeta$. It can be shown:

**Proposition 6.** *There is a recursive operator $\Gamma$ with $\xi(\Gamma) = \zeta$.*

**Proof:** The universal Infinite Time Turing machine of Hamkins and Kidder ([6]) is in effect such a recursive operator $\Gamma$. That $\xi(\Gamma) = \zeta$ is shown in [16]. Q.E.D.
   Some quasi-inductive definitions may reach a fixed point.

**Definition 7.** *We say that $\Gamma$ reaches a fixed point on $X$, if there is $\alpha$ so that $\Gamma_\alpha(X) = \Gamma_{\alpha+1}(X)$; and if so we call $\Gamma_\alpha(X)$ the fixed point.*

**Proposition 8.** *For any arithmetical operator $\Gamma$, either $\xi(\Gamma) = \zeta$, or else there is an equivalent recursive operator $\tilde{\Gamma}$ which reaches $\Gamma_{\xi(\Gamma)}(\varnothing)$ as a fixed point, that is there is a recursive operator $\tilde{\Gamma}$ with: $\exists \alpha < \zeta \ \tilde{\Gamma}_\alpha(\varnothing) = \tilde{\Gamma}_{\alpha+1}(\varnothing) = \Gamma_{\xi(\Gamma)}(\varnothing).$*

The quasi-inductive definitions that reach fixed points (on $\varnothing$, or on some particular input $X$) form an interesting subclass. Investigation of such is an appealing combination of admissibility theory and reflection properties of ordinals.

Propositions 6 and ? indicate that in one sense, to study recursive operators, is to study all arithmetic ones: if one has a $\Sigma_n$-definable operator, one seemingly only needs to look instead at $\Pi_{n+1}$-reflecting ordinals. For such an operator $\Gamma$ as in Proposition 4, it is easy to see from the definition of $(\zeta, \Sigma)$ that $\Gamma_\zeta(\varnothing) = \Gamma_\Sigma(\varnothing)$ and thus we might call $(\zeta, \Sigma)$ a ''repeat pair'' for $\Gamma$ and $\varnothing$. Again one may show for such a $\Gamma$, that $(\zeta, \Sigma)$ is the lexicographic least such repeat pair. We use this to formulate a definition allowing parameters $x$ as starting inputs.

**Definition 9. AQI** *is the sentence: ''For every arithmetic operator $\Gamma$, for every $x \subseteq \omega$, there is a wellordering $W$ with a repeat pair $(\zeta(\Gamma, x), \Sigma(\Gamma, x))$ in* $\mathrm{Field}(W)$''*. If an arithmetic operator $\Gamma$ acting on $x$ has a repeat pair (in* $\mathrm{Field}(W)$)*, we say that $\Gamma$ stabilizes along $W$ (with input $x$).*

Clearly a certain amount of set theory (or analysis) is needed to show that every operator stabilizes. Reformulated using Prop. 3, this is thus:

**Lemma 10.** $\mathrm{KP} \vdash \mathbf{AQI} \quad \longleftrightarrow \quad \forall x \subseteq \omega \ \exists \xi, \sigma \, ( \ L_\xi[x] \prec_{\Sigma_2} L_\sigma[x] \, ).$

We note some facts concerning the pair $(\zeta, \Sigma)$ in $L$:

**Proposition 11.** *(i)* ([16] *Thm. 2.1) $L_\zeta$ is a model of $\Sigma_2$-KP (and is a union of such).*
*(ii)* ([17] *Cor.3.4) $L_\Sigma$ is a model of $\mathrm{KPI}_0 + \Sigma_2$-Separation$_0$, but not of* KP.

In the language of subsystems of second order number theory, the reals of $L_\Sigma$ form a $\beta$-model of $\mathbf{\Pi_1^1}$-$\mathbf{CA_0}$ together with a $\beta_3$-submodel (provided by the reals of $L_\zeta$), and in fact it is the minimum such $\beta$-model. (See [14] VII.7.1 for the notion of a $\beta = \beta_1$ and $\beta_k$-submodels.)

Prop. 9 (i) then already shows that **AQI** is stronger than $\mathbf{\Delta_3^1}$-$\mathbf{CA_0}$ since any $\Sigma_2$-extendible is already a (union of) models of $\Sigma_2$-KP. Somewhat more formally:

**Corollary 12.** $\mathbf{\Delta_3^1\text{-}CA_0}+\mathbf{AQI}\vdash \forall x(\exists \beta\text{-model of } \mathbf{\Delta_3^1\text{-}CA_0} \text{ containing } x)$.

Proof: In fact we have: $\mathbf{\Delta_3^1\text{-}CA_0}+\mathbf{AQI}\vdash$ "*There is a countably coded $\beta$-model of $\mathbf{\Delta_3^1\text{-}CA_0}$*". This is because the former proves "$\forall x \exists W \in$ WO$(\exists z, s \in \text{Field}(W) \wedge L_{\|z\|}[x] \prec_{\Sigma_2} L_{\|s\|}[x])$." By Prop. 9 (i), $\mathcal{P}(\mathbb{N})$ of the model $L_{\|z\|}[x]$ form a model of $\mathbf{\Delta_3^1\text{-}CA_0}$.    Q.E.D.

In Sections 2,3 we prove the principle Theorems 2, and 1 of the abstract respectively. These with the last corollary discharge all the abstract's assertions. (The reader may consult [14] on the mutual interpretability between theories in subsystems of second order number theory, such as $\Delta_3^1$-CA$_0$ for example, and set theoretical counterparts $\Delta_3^1$-CA$_0^{\text{set}}$, or $\Sigma_2$-KP. We do not have to do any fine analysis of provability in any such subsystems, as we are typically showing that a rather strongtheory $A$ proves the existence of (many) $\beta$-models of theory $B$.)

## 1.3  Strategies and game trees.

We assume familiarity with the basic notions of two person perfect information games played using integers. We shall follow Martin and shall assume that such games are played on *game trees* $T \subseteq {}^{<\omega}\omega$ although we disallow terminal nodes. We let $G(A;T)$ denote the game with payoff set $A \cap \lceil T \rceil$ where $\lceil T \rceil$ denotes the set of all plays in $T$. A *position* in the game is simply a finite sequence $r \in T$. For $q \in T$ we let $T_q$ denote the set of all positions $r \in T$ where $r \supseteq q$.

**Lemma 13.** *Let $A$ be arithmetic; let $M$ be a transitive model of* KPI$_0$ *with $T \in M$. Then (i) "$G(A;T)$ is not a win for $I$" is $\Pi_1^M$; (ii) if this holds then "$p$ is a position in $II$'s non-losing quasi-strategy for $G(A;T)$" is $\Pi_1^M$.*

**Proof.** "$G(A;T)$ is not a win for $I$" is equivalent to "$\forall \sigma \in {}^\omega\omega(\sigma$ is not a winning strategy for $I$ in $G(A;T)$ )"; which in turn is equivalent to : "$\forall \sigma($if $\sigma$ is a strategy for $I$ in $G(A;T)$ then $\exists r \in {}^\omega\omega \,\sigma * r \notin A \cap \lceil T \rceil)$". The set $\{r \mid \sigma * r \notin A \cap \lceil T \rceil\}$ is then $\Delta_1^1(\sigma, T)$, and hence, if non-empty with $\sigma \in M$, has an element in $M$. This completes (i). For (ii): let $T'$ be $II$'s non-losing quasi-strategy for $G(A;T)$. Then $p \in T' \Leftrightarrow \text{lh}(p) = k \to \forall n \le k \,(p \upharpoonright n \in T) \wedge \forall 2n+1 < k(q = (p_0, p_1, ..., p_{2n+1}) \to$ "$G(A;T_q)$ is not a win for $I$."    $\square$

One should note that a quasi-strategy for $II$ in $G(A;T)$ is then a subtree of $T$ that does not restrict $I$'s moves in any way, only $II$'s moves.

# 2 AQI is weaker than $\Sigma_3^0$-Determinacy.

**Definition 14.** $T_\alpha^j(X)$ *denotes the set of* $\Sigma_j$ *formulae, true of parameters from* $X$, *in the structure* $\langle J_\alpha, \in \rangle$; $T_\alpha^j$ *abbreviates* $T_\alpha^j(\varnothing)$, *the* $\Sigma_j$-*theory of* $\langle J_\alpha, \in \rangle$.

**Definition 15.** *Let* $E_0$ *be the class of* $\Sigma_2$-*extendible ordinals. If* $E_k$ *is defined, let* $E_{k+1}$ *be the class* $E_0 \cap E_k^*$. *Let* $E = \cap_{k \in \omega} E_k$.

The classes $E_k$ we can think of as having depth $k$ in the ''$\Sigma_2$-extendible limits of $\Sigma_2$-extendible ...'' hierarchy: if $\gamma \in E_k$ then there are ordinals $\gamma = \mu_k \leq \mu_{k-1} \leq ... \leq \mu_0 < \nu_0 < \nu_1 < ... < \nu_k$ satisfying $L_{\mu_j} \prec_{\Sigma_2} L_{\nu_j}$ for $j \leq k$.

**Theorem 16.** $\Sigma_2$- KP $+ \forall \alpha \exists \beta, \gamma (\alpha < \beta < \gamma \wedge L_\beta \prec_{\Sigma_2} L_\gamma) \not\vdash \Sigma_3^0$-Determinacy.

**Proof.** We shall show that the least level of the $L$-hierarchy that is a model of the antecedent theory is a model $\mathcal{M}_0 = L_\delta$ in which $\Sigma_3^0$-Determinacy fails. The reals of this model form a $\beta$-model of $\mathbf{\Delta_3^1}$-$\mathbf{CA_0} + \mathbf{AQI}$. The rest of this section is taken up with proving this theorem.

In fact we shall prove something slightly stronger in the form of the following Lemma:

**Lemma 17.** $\Sigma_2$- KP $+ \forall \alpha \exists \beta (\alpha < \beta \wedge \beta \in E^*) \not\vdash \Sigma_3^0$-Determinacy.

We do this, using a technique that goes back to H.Friedman, by defining certain games $G_\psi$ so that codes for initial segments of the $L$-hierarchy are recursive in any winning strategy for the game. So henceforth, let $\mathcal{M} = L_\delta$ be the least level of the antecedent theory in the statement of the lemma.

Let $\Psi = \{ \psi \mid \psi \in \Sigma_1 \wedge L_\delta \vDash \psi \} = T_\delta^1$ be the $\Sigma_1$-theory of $L_\delta$.[2]

(1) ''$\Sigma_3^0$-Determinacy'' is $\Sigma_1^{\mathrm{KPI}}$.

**Proof.** $\Sigma_3^0$-Determinacy is equivalent to

---

$\forall n \in \omega$ [if $A_n$ is the $n$'th $\Sigma_3^0$ set then $\exists \sigma$ ($\sigma$ is a winning strategy for a player in $G(A_n; {}^{<\omega}\omega))$]

The statement that $\sigma$ is a winning strategy for Player $I$ is equivalent to saying: "The $\Delta_1^1(\sigma)$ set $\{x \mid \sigma * x \notin A\}$ is empty". If it were non-empty then it would have a member $\Sigma_1^1(\sigma)$-definable and thus definable over the least admissible set containing $\sigma$. Thus "$\sigma$ is a winning strategy" is thus $\Sigma_1^{\mathrm{KPI}}$. We thus have a numerical quantifier in front of a $\Sigma_1^{\mathrm{KPI}}$ predicate, and is thus overall $\Sigma_1^{\mathrm{KPI}}$.                                  $\square$

Hence, were $\Sigma_3^0$-Determinacy to hold in $L_\delta$ it would be equivalent to some $\psi \in \Psi$. For any $\psi \in \Psi$ we define: $\alpha_\psi =$ the least $\beta$ so that $L_\beta \vDash \mathrm{KP} + \psi$.

**Note 18.** The leastness of $\beta$ ensures that every $x \in L_{\alpha_\psi}$ is $\Sigma_1$-definable by some parameter free $\Sigma_1$ term $t_x$. (In other words the $\Sigma_1$-skolem hull inside $\langle L_{\alpha_\psi}, \in \rangle$ of $\varnothing$ is all of $L_{\alpha_\psi}$ itself.)

The following is straightforward:

(2) Let $\bar{\alpha} = \sup \{\alpha_\psi \mid \psi \in \Psi\}$. Let $\alpha' =$ the least $\beta(L_\beta \prec_{\Sigma_1} L_\delta)$. Then $\alpha' = \bar{\alpha}$.

We shall show for every $\psi \in \Psi$ there is a game $G_\psi$ with a $\Pi_3^0$ payoff set, but without a winning strategy in $L_{\alpha_\psi}$. In view of the comment just before (2), this will suffice.

Fix for the rest of the argument $\psi \in \Psi$. Let $\alpha$ denote $\alpha_\psi$. We consider the following game $G = G_\psi$.

| $I$ | plays | $m_0, \ m_1, \ldots, m_i$ | $x = (m_0, m_1, \ldots, m_i, \ldots)$ |
|-----|-------|---------------------------|----------------------------------------|
| $II$ | plays | $n_0, \ n_1, \ldots, \ n_i$ | $y = (n_0, n_1, \ldots, n_i, \ldots)$ |

in the usual way, playing in the $i$'th round integers $(m_i, n_i)$. Let $z = x \oplus y$.

*Rules for I.*

Let $T$ be the theory $\mathrm{KP} + V = L + \psi$. $x$ must be a set of Gödel numbers for the complete $\Sigma_1$-theory of an $\omega$-model of $T +$ "*there is no set model of $T$*".

Using the Note 18 we denote by $\langle M, E \rangle$ the model $I$ essentially constructs if he obeys this rule. We may regard also as part of the rule that $x$ as given by $I$ should be specified simply by $I$ stating ``$k \in T_M^1$'' or ``$k \notin T_M^1$'' where $T_M = T_M^1$ is the standard $\Sigma_1$-code for the appropriate level of the $L$-hierarchy. Also, as in $L_{\alpha_\psi}$, in $\langle M, E \rangle$, every set is given by a $\Sigma_1$ parameter free skolem term.

**Note 19.** If $\langle M, E \rangle$ is wellfounded then it is isomorphic to $\langle L_{\alpha_\psi}, \in \rangle$.

Amongst the codes for sentences that $I$ plays are those of the form

$$\ulcorner t_m \in \mathrm{On} \wedge t_n \in \mathrm{On} \wedge t_m < t_n \urcorner$$

These we shall use to formulate rules for player $II$. So far the Rules for $I$ amount to a $\Pi_2^0$ condition on $x$ and so on $z$. (We may take a recursive listing of $\Sigma_1$-sentences $\langle \psi_k | k \in \omega \rangle$ and we then require $\forall k \exists k'(m_{k'} = \ulcorner \psi_k \urcorner \vee m_{k'} = \ulcorner \neg \psi_k \urcorner)$, thus the theory $I$ constructs will be $\Sigma_1$-complete; we obtain that $M$ has at least the integers as standard also by a $\Pi_2^0$ condition.) Let $r : \omega \to \omega \times \omega$ be a recursive enumeration of $\omega^2$ in which each $(i, j)$ appears infinitely often.

*Rules for II.*

At round $k$:

if $(i, j) = r(k)$ and $n_k \neq 0$, then we shall say that ``*II makes the entry $n_k$ on list $L_{i,j}$*.'' These 'Listing' *Rules* here require her to list terms in a correct order. She may make an entry on list $L_{i,j}$ in round $k$ if:

*Either* $L_{i,j}$ is empty at the current round, in which case $n_k$ can be any term $t_s$ *as long as I has asserted at an earlier round* $\ulcorner t_s \in \mathrm{On} \urcorner \in T_M$;

*or* $L_{i,j} \neq \varnothing$, and if $t_s$ was the last entry $II$ made on this list, then $n_k$ can be any term $t_r$, again *provided that I has at an earlier round* $k' \leq k$ asserted $m_{k'} = \ulcorner t_r \in \mathrm{On} \wedge t_r < t_s \urcorner \in T_M$.

*The winning conditions.* $I$ wins immediately at a finite round if $II$ breaks one of her *Listing Rules* just enumerated. $II$ wins if $I$ fails to obey his conditions on $x$, or both players obey their respective rules and additionally

$\exists (i, j)$ [*II makes infinitely many entries on list $L_{i,j}$*].

This is a $\Sigma_3^0$ winning condition for $II$ on $z$. Hence $G_\psi$ has a $\Pi_3^0$ payoff set.

**Remark 20.** In other words, if $I$ obeys his rules, $II$ can win if for some $(i, j)$, $r^{-1}``(i, j)$ in effect picks out an infinite descending chain through the ordinals of the model $\mathcal{M}$ that $I$ reveals *via* the gödel numbers of the $\Sigma_1$ sentences true in $\mathcal{M}$.

**Remark 21.** $II$ is not allowed to make an entry indicating that $t_s < t_r$ until $I$ has asserted this at some earlier stage. $II$ is thus not predicting what the model will look like below $t_r$; by making an entry on a list she is merely adverting to the fact that $I$ has revealed that $t_s < t_r$.

**Lemma 22.** *$I$ has a winning strategy.*

**Proof:** $I$ plays out all ``$k \in T_M$'' for all $k \in T_{\alpha_\psi}^1$, and ``$k \notin T_M$'' for all $k \notin T_\alpha^1$. Obviously then, $\langle M, E \rangle \simeq \langle L_{\alpha_\psi}, \in \rangle$ and $II$ has no chance to pick out any infinite descending chains.  Q.E.D. Lemma

The point is the following:

**Lemma 23.** *Let $\tau$ be any winning strategy for $I$. Let $x = T_{\alpha_\psi}^1$; then $x \leq_T \tau$.*

From this the theorem then follows as $x \notin L_{\alpha_\psi}$,  being essentially the latter's $\Sigma_1$-truth set.

**Proof of Lemma 23** We argue that, with $II$ only playing constantly $n_k = 0$ for all $k$, that $I$ is forced to play for $x$ a list of all the correct facts ``$k \in / \notin T$'' for $T_M = T_{\alpha_\psi}^1$. The point is to show that if at any time $I$ deviates from this course of action, then he will lose - and hence the purported strategy $\tau$ is not a winning one.

$II$ plays ``*Pass*'' ($n_k = 0$) until such a point, if ever, when $I$ asserts $k \in T_M$ or $k \notin T_M$ whereas in reality $k \notin T_{\alpha_\psi}^1$ or $k \in T_{\alpha_\psi}^1$. At this point $II$ knows that $I$'s eventual model $\langle M, E \rangle$ will be illfounded, and so she must act to discover a descending chain. In this case we shall set $\beta = \beta_M =_{df}$ On $\cap$ WFP$(M)$. However she will not yet know, and in fact will not at any move know where $\beta_M$ lies. As $(KP)_M$ by the Truncation Lemma (*cf.* [1]) $\beta_M \in$ ADM. By our requirements on the theory $T_{\alpha_\psi}^1$, and upwards persistence of $\Sigma_1$ formulae, we must have $\beta_M \leq \alpha_\psi$.

**Definition 24.** *Let* $F: \omega \twoheadrightarrow \text{ADM} \cap \alpha_\psi + 1$ *be some fixed surjection.*

The idea is that at rounds $k$ where $r(k) = (i, j)$ *II* will be making the working assumption that the ordinal height of the wellfounded part of $M$, $\beta_M$, is precisely $F(i)$, and will be trying to find an illfounded chain through $\text{On}^M$ above $\beta_M$. She will be working simultaneously on all such possible $\beta_M$. We shall prove that if $I$ deviates from enumerating $T_{\alpha_\psi}^1$, knowing that one of them is the correct assumption, she can nevertheless be successful and win the game $G_\psi$; thus $I$ is forced to play only the truth concerning $T_{\alpha_\psi}^1$.

We assume then that $I$ has played an untruth. We concentrate on a fixed $i$ and hence on $\beta = \beta_M = F(i)$, and describe how *II* can move in rounds $k$ with $r(k) = (i, j)$.

(3) *Claim* $\exists \bar{a} \notin \text{WFP}(M) \forall b < \bar{a} \, (b \notin \text{WFP}(M) \rightarrow T_b^2 \not\subset \tilde{T} = T_\beta^2$.

**Proof.** Supposed this failed, then $\forall \bar{a} \notin \text{WFP}(M) \exists b < \bar{a} \, (b \notin \text{WFP}(M) \wedge T_b^2 \subset \tilde{T}$. For such $b$ we shall show that $J_{\gamma_b} \prec_{\Sigma_2} J_b$ for a $\gamma_b < \beta$. To do this we could employ a version of Theorem 4 on Uniformly Definable Skolem functions, to apply to the non-wellfounded model $M$. To do this one would simply apply $h_b^2$ (the uniformly defined $\Sigma_2$-skolem function, as defined over the illfounded model $J_b$), and look at $h_b^2``\omega \times \omega$, and argue that this is a $\Sigma_2$ skolem hull of $J_b$ which is transtive, and in fact wellfounded, and is then a $J_{\gamma_b}$ for a $\gamma_b < \beta$. However we prefer to argue for this directly as follows.

Let $\bar{\eta} =_{\text{df}} \sup \{c < b \mid \exists f \in \Sigma_2^{J_b}, f: \omega \to c, f \text{ partial, onto}\}$. We first claim that $\bar{\eta} < \beta$. Clearly equality fails, as otherwise that would make $\beta$ definable inside $M$ from $b$. If however $c \notin \text{WFP}(M)$, with $f \in \Sigma_2^{L_b}, f$ partial, but onto $c \not\leq \beta$, then the sentences ``$f(n)\downarrow, f(m)\downarrow \wedge f(n) < f(m) \in$ On'' are all in $T_b^2$ and so in $\tilde{T}$. This is absurd as $\beta$ is wellfounded!

As there is a $\Delta_1^{J_{\bar{\eta}}}$ map of $\bar{\eta}$ onto $J_{\bar{\eta}}$, we could equally well write: $J_{\bar{\eta}} = \bigcup \{J_\alpha \mid \exists f \in \Sigma_2^{L_b}, f: \omega \to J_\alpha, f \text{ partial, onto}\}$. We wish to claim that $J_{\bar{\eta}} \prec_{\Sigma_2} J_b$. Suppose $J_b \vDash \exists u \psi(u, \vec{\xi})$ where $\vec{\xi} < \bar{\eta}$ and $\psi \in \Pi_1$.

As $\vec{\xi} = f_0(n)$ for some $\Sigma_2^{J_b} f_0$, we have $J_b \vDash \exists y [\exists y (y = f_0(n) \wedge \psi(u, y)]$. This is a $\Sigma_2$ sentence and is in $T_b^2$. Suppose $y = f_0(n) \leftrightarrow \exists v \varphi(v, y, n)$ say for some $\Pi_1$ $\varphi$. If $\sigma \in T_{\bar{\eta}}^2$ and ``$\vec{\xi} = f_0(n)$'' holds in $L_{\bar{\eta}}$ then we are done. Let $\delta$ be the least ordinal such that $\forall \delta' > \delta L_{\delta'} \vDash \sigma$. Then $\{\delta\} \in \Pi_1^{L_b}$, and thus $\delta < \bar{\eta}$. Hence $\sigma \in T_{\bar{\eta}}^2$. We now consider the $\Pi_1$ $\varphi$. Let $\gamma$ be least so that for some $m < \omega$, $\varphi((h^1(m, \gamma))_0, (h^1(m, \gamma))_1, n) \wedge J_\gamma \vDash$ ''$\forall m \forall \gamma' \neg \exists z(z = h^1(m, \gamma') \wedge \varphi((z)_0, (z)_1, n))$.''

This is a $\Pi_1^{L_b}$ expression about $\gamma$. Now note that for the least such $\gamma$ satisfying the above, we must have $J_\gamma \prec_{\Sigma_1} J_b$. (For otherwise $\gamma$ itself is of the form $h^1(m', \gamma')$ for a $\gamma' < \gamma$ and the lesser $\gamma'$ can be used in place of $\gamma$.) We further note then by this $\Sigma_1$-reflection on the quoted formula in the final conjunct, that $J_b$ itself is a model of ''$\forall m \forall \gamma' < \gamma \neg \exists z(z = h^1(m, \gamma') \wedge \varphi((z)_0, (z)_1, n))$.'' The conclusion is that the least such $\gamma$ so that ''$\exists z(z = h^1(m, \gamma) \wedge \varphi((z)_0, (z)_1, n))$'' is expressible in a $\Pi_1^{J_b}$ way. But that means $\gamma < \bar{\eta}$. However then $J_{\bar{\eta}}$ is a model of this statement and so ''$\vec{\xi} = f_0(n)$'' holds in $J_{\bar{\eta}}$.

Hence for such a $b$ we have $(J_{\gamma_b} \prec_{\Sigma_2} J_b)_M$. However the supposition implies there is an infinite descending chain of such $b$ in the illfounded part of $M$. This implies that we have an infinite nested sequence of $\Sigma_2$ intervals: there exists $\langle b_n | n < \omega \rangle$, $\langle \gamma_n | n < \omega \rangle$ with $(\gamma_n \leq \gamma_{n+1} \leq ... < b_{n+1} < b_n)$, and with $(J_{\gamma_n} \prec_{\Sigma_2} J_{b_n})_M$, for $n < \omega$. This implies that each $\gamma_n \in E$, and in fact in $E^*, (E^*)^*, ...$thus contradicting our smallness hypothesis.     Q.E.D. (3)

Let $\langle t_k | k \in \omega \rangle$ be our priorly fixed recursive enumeration of the $\Sigma_1$-skolem terms (using the standard $\Sigma_1$-skolem function, this could simply be an enumeration of $\langle h^1(i, n) \mid i, n < \omega \rangle$). $II$ makes the additional working assumption, or guess if you will, that $t_j^M = a_0$, where $a_0$ is a witness for $\bar{a}$ to the truth of the last Claim. (Again the point is that $II$ does not know in advance which term in $M$ will denote such $a_0$.) As $I$ reveals more and more facts about his model, he must, if $M$ is not going to be isomorphic to $L_{\alpha_\psi}$, at some point reveal a $\Sigma_1$-fact which is true in $M$ but false in $L_{\alpha_\psi}$. There really is then such an $M$-ordinal $a_0$. $II$ will, in effect, place her 'guess' $a_0 = t_j^M$ at the head of her putative descending chain, on list $L_{i,j}$. In order to choose the next element of the chain $II$ considers the set $\tilde{T} = T_{\bar{\beta}}^2$. Set $T_0 = (T_{t_j}^2)_M$.

$II$ now waits until $I$ asserts that some $\sigma_0$ is in $T_0$, (this itself being one of the $\Sigma_1$ facts about $M$ she must enumerate) but $II$ sees is not in $\tilde{T}$. (If $II$ is wrong in her guess about $t_j$ of course, then she may fruitlessly wait for ever...)

(4) Suppose $M \vDash$ " $a_1 < a_0$ is least so that $\forall b \le a_0 (b \ge a_1 \rightarrow (\sigma_0)_{L_b})$." Then $a_1 \notin \text{WFP}(M)$

**Proof**: Were $a_1 \in L_\beta$ then we should have $\sigma_0 \in \tilde{T}$.      Q.E.D. (4)

*II* may thus wait until *I* asserts that some such $\sigma \in T_0 \setminus \tilde{T}$ and additionally, perhaps later, the $\Sigma_1$ fact that some term $t_{j_1}$ names the ordinal $a_1$ defined in (4) above. And at some round $l$ then, *I* must play the number $m_l = \ulcorner t_{j_1} \in \text{On} \wedge t_j \in \text{On} \wedge t_{j_1} < t_j \urcorner$; once all these facts have been gathered together, Agathe may at the next appropriate round $k$ with $r(k) = (i, j)$, set $n_k = t_{j_1}$.

*II* now has two elements of a descending chain in the illfounded part of $M$. Now she watches out for assertions that *I* makes about $T_1 = (T_2^{t_{j_1}})_M$, waiting for some $\sigma_1$ asserted by him to be in $T_1$ but which does not lie in $\tilde{T}$. By exactly the same considerations that held at (4) some $a_2, t_{j_2}$, are definable, and so she can continue. By the end of the game, *if* this working assumption about $\beta_M$ and $t_j$ was the correct one, the chain so defined by continuation of this process will be infinite, and she will have won.

If *I* deviates from playing the correct truth set, then at least one of *II*'s assumptions will turn out to be a correct one, and hence she will be assured of winning.      Q.E.D. **(Lemma 23 & Theorem 16)**                   □

**Remark 25.** One might wonder why $\Sigma_2$ theories play the role they do here. However in this kind of game they are necessary. For suppose we tried to do without considering $\Sigma_2$ sentences but *II* tried to divine a descending chain using only $\Sigma_1$ sentences. Any $M$-ordinal is defined by a $\Sigma_1$ term $t_k$ say. However 'new' $\Sigma_1$ facts about standard ordinals $< \beta_M$ may be true in $M$. Hence *II* cannot use the terms alone to pick out a chain without the danger of falling into the wellfounded part. Even if she waits to hear that a term $t_l$ specifies a level hopefully in the ill-founded partwhere the $\Sigma_1$ term $t_k$ (*i.e.* one not mentioned in $\tilde{T}$) becomes true (which may be at a much higher level than the ordinal named by $t_k$) then this could also fail. For suppose that in $M$ there is a $\gamma < b$ with $\gamma$ in the wellfounded part, and $b$ in the illfounded part, and with $J_\gamma \prec_{\Sigma_1} J_b$. No new $\Sigma_1$ facts become true within the interval $[\gamma, b)_M$. *I* may artfully play terms so that *II* ends up choosing a finite chain of terms naming ordinals where new $\Sigma_1$ facts are witnesed, with bottom element $b$. If she drops any further using such terms, she will be in the wellfounded part. Hence the use of the device above of looking where $\Sigma_2$-sentences become true.

**Remark 26.** One might further wonder whether the use of infinitely many lists only indicates our poor ability to devise a better game that used only one list. However the use here of infinitely many lists is necessary: for suppose it could be effected with a single list say, then the game would be $\Pi_2^0 \backslash \Pi_2^0$ (at least in the version described here) and it is known (see [7] that the least level $L_\delta$ with $L_\delta \vDash \Sigma_1$-Separation is a model of Boolean($\Sigma_2^0$)-Determinacy, and this is a very much smaller ordinal than the first $\Sigma_2$-extendible $\zeta$ above.

Recall the following definition.

**Definition 27.** *Let $\Gamma$ be a pointclass. A set $Z \subseteq \mathbb{N}$ is said to be in $\partial\Gamma$ if there is a set $X \subseteq \mathbb{N} \times \mathbb{N}^\mathbb{N}$ in $\Gamma$ so that*
$$Z = \{n \,|\, I \text{ has a winning strategy in } G(X_n; \omega^{<\omega})$$
*where $X_n = \{y \,|\, \langle n, y \rangle \in X\}$.*

The next Corollary is in a sense merely a restatement of the result above.

**Corollary 28.** $\bar{\alpha} \le \sigma_3$. *Indeed each $T^1_{\alpha_\psi}$ is a $\partial\Sigma_3^0$ set of integers.*

**Proof:** Let $\alpha_\psi$ *etc.* be defined as above. We switch roles in the games. $I$ will try and find descending chains through $II$'s model $M$. (This is only to make the payoff set $\Sigma_3^0$ rather than $\Pi_3^0$.) For $\varphi \in \Sigma_1$ let $G_{\psi,\varphi}$ be the game described in the last theorem, except that $II$ must now play a code $x$ for a model of $T$ + ''there is no set model of $T$'' $+\neg\varphi$. Everything else remains the same *mutatis mutandis: I*'s task is to find an infinite descending chain through the ordinals of $II$'s model. Note that if $\varphi \in T^1_{\alpha_\psi}$ $I$ now has a winning strategy: for if $II$ obeys her rules, and $x$ codes an $\omega$-model $M$ of this theory, then $M$ is not wellfounded, and has WFP$(M) \cap$ On $< \alpha_\varphi$, where the latter is the least admissible $\alpha$ where $\varphi$ is true in $L_\alpha$. However $I$ playing (just as $II$ did in the last theorem) can find a descending chain and win. On the other hand if $\varphi \notin T^1_{\alpha_\psi}$ $II$ may just play a code for the true wellfounded $L_{\alpha_\psi}$ and so win. This shows that $T^1_{\alpha_\psi}$ is a $\partial\Sigma_3^0$ set of integers.

Suppose now that $\bar{\alpha} > \sigma_3$. Let $\psi$ be such that $\alpha_\psi$ is the second least admissible greater than $\sigma_3$. There is thus a set $H \in L_{\alpha_\psi}$ (definable over the next admissible $\gamma > \sigma_3$) containing winning strategies for all $\Sigma_3^0$-games that are a win for player $I$ and in particular a set $H_0$ of those winning strategies for $I$ in games of the form $G_{\psi,\varphi}$. Hence membership of $\varphi$ in $T^1_{\alpha_\psi}$ is determined by searching through $H_0$ for a winning strategy for $I$; this is a $\Sigma_1$-search. Hence $T^1_{\alpha_\psi} \in \Delta_1^{L_{\alpha_\psi}}(\{H_0\})$. Hence $T^1_{\alpha_\psi} \in L_{\alpha_\psi}$ which is a contradiction.          Q.E.D.

**Corollary 29.** *The complete AQI set is a $\eth \Sigma_3^0$ (but not $\eth \Sigma_2^0$) set of integers.*

# 3 (Boldface) $\mathbf{\Sigma_3^0}$-Determinacy is weaker than $\Sigma_2$-Separation.

We shall closely follow Martin's account of Davis' proof ([3]) of $\Sigma_3^0$-Determinacy. That account is performed within $\mathrm{ZC}^- + \Sigma_1$-Replacement, but we shall pay attention to definability considerations. It is fairly easy to see that the proof can be effected with $\Sigma_2$-Replacement and $\Sigma_2$-Separation, but we want to do better than that. As remarked above $\mathrm{KP} + \Sigma_2$-Separation proves the existence of models $M$ of $V = L +$ *"There are ordinals $\gamma_0 < \gamma_1 \in \mathrm{ON}^M$ with $(L_{\gamma_0} \prec_{\Sigma_1} V \wedge L_{\gamma_0} \prec_{\Sigma_2} L_{\gamma_1})$"* $+ \forall x \exists y (\mathrm{Trans}(y) \wedge y \text{ is admissible})$.

Notice that if $M$ is a model as described, then for parameters $t \in L_{\gamma_0}$, $\Sigma_2^{L_{\gamma_1}}(\{t\})$-definable subsets of $\omega$ are all in $L_{\gamma_1}$, and are in fact so-definable over $(L_{\gamma_0})_M$. Suppose we set $\gamma_2 = \mathrm{ON}^M$. We shall show for such an $M$:

**Theorem 30.** $L_{\gamma_0}^M \vDash \mathbf{\Sigma_3^0}$-Determinacy.

By "$N \vDash \mathbf{\Sigma_3^0}$-Determinacy" we shall mean that for game trees $T \subseteq {}^{<\omega}\omega$ with $T \in N$ and for any $A \in \Sigma_3^0(T)$, that a winning strategy $\sigma$ for player $I$ or $II$ exists in $N$. As $L_{\gamma_0}, L_{\gamma_1}$ are both models of $\mathrm{KPI}_0$, the property of any $\sigma \in L_{\gamma_0}$ being a winning strategy for such a game $G(A; T)$ is absolute between $L_{\gamma_0}$ and the universe. This will complete the first theorem of the abstract (as well as the non-reversibility of its implication, since by taking $\gamma_0$, $\gamma_1, \gamma_2$ least with such properties we have that then $\Sigma_2$-Separation fails in $L_{\gamma_0}$).

**Proof:** We shall assume that $V = M$ where $M$ is a model with the above properties. We shall thus drop the subscript $M$ throughout the proof. Let $A$ be $\Sigma_3^0(T)$ for some game tree $T \in L_{\gamma_0}$. We shall show $G(A; T)$ is determined, with either a winning strategy for $II$ in $L_{\gamma_0}$ or a winning strategy for $I$ definable over - somewhat loosely speaking- $L_{\gamma_2}$. We suppose that neither player has a winning strategy in $L_{\gamma_0}$ and shall prove that $I$, does. The key lemma is the following.

**Lemma 31.** *Let $B \subseteq A \subseteq \lceil T \rceil$ with $B \in \Pi_2^0$. If $(G(A; T)$ is not a win for $I)_{L_{\gamma_0}}$, then there is a quasi-strategy $T^* \in L_{\gamma_0}$ for $II$ with the following properties:*
*(i) $\lceil T^* \rceil \cap B = \varnothing$ ;*
*(ii) $(G(A; T^*)$ is not a win for $I)_{L_{\gamma_0}}$.*

**Remark**: We recall from Lemma 13 that "$G(A;T)$ is not a win for $I$" is a $\Pi_1^{\mathrm{KPI}_0}$ statement about $A$, $T$ and is absolute between $L_{\gamma_0}$, $L_{\gamma_1}$, and $V = L_{\gamma_2}$ by our elementarity properties assumed of these models; similarly if there is such a $T^*$ as claimed by the conclusion, then $(G(A;T^*)$ is not a win for $I)_{L_{\gamma_i}}$ for $i < 3$.

**Proof of Lemma.** If $T'$ is $II$'s non-losing quasi-strategy for $G(A;T)$, as defined over $L_{\gamma_0}$, then by Lemma 13(ii) membership in "$p \in T'$" is not only $\Pi_1^{L_{\gamma_0}}$, but also $\Pi_1^{L_{\gamma_i}}$ for $i = 1, 2$, due to the assumed $\Sigma_1$ elementarity. As $L_{\gamma_0}$ is a model of $\Sigma_1$-Comprehension, $T' \in L_{\gamma_0}$, and we thus have that for every $p \in T$ :

(1) $(p \in T')_{L_{\gamma_0}} \longleftrightarrow (p \in T')_{L_{\gamma_i}}$ *for $i = 1, 2$. Hence $T'$ is $II$'s non-losing quasi-strategy for $G(A;T)$ in $V$.*

Likewise for any $p \in T'$ "$G(A, T'_p)$ is not a win for $I$" is absolute between $L_{\gamma_i}$ for $i < 3$.

Following closely the original argument, we call a position $p \in T'$ *good* if there is a quasi-strategy $T^*$ for $II$ in $T'_p$ so that the following hold:

(i) $\lceil T^* \rceil \cap B = \varnothing$;
(ii) $G(A;T^*)$ is not a win for $I$.

Here (ii), again, is a $\Pi_1^{\mathrm{KPI}_0}(T^*)$ statement by Lemma 13. Thus the existence of such a $T^*$ becomes a $\Sigma_2$ fact about $T$ (as (i) is simply $\Delta_1$). Hence by our assumed $\Sigma_2$ reflection properties it makes no difference whether we defined 'good' relative to $L_{\gamma_0}$ or $L_{\gamma_1}$:

(2) "$p$ is good" is $\Sigma_2^{\mathrm{KPI}_0}(T)$ and hence $(p \text{ is good})_{L_{\gamma_0}} \longleftrightarrow (p \text{ is good})_{L_{\gamma_1}}$

We set $H \subseteq {}^{<\omega}\omega$ to be the class $\{p \in T' \mid (p \text{ is good})_{L_{\gamma_0}}\}$. Then

(3) $H \in \Sigma_2^{L_{\gamma_0}}(T)$ and hence $H$ is a set in $L_{\gamma_1}$.

We define the function $t \colon H \longrightarrow L_{\gamma_0}$ by:

$t(p) = L$-least quasistrategy $\hat{T}(p)$ witnessing (i) and (ii) that $p$ is good.

Then $t$ is definable over $L_{\gamma_0}$ in $T'$, but is also $\Sigma_2^{L_{\gamma_1}}(\{T'\})$ as a relation, and, noting that $H$, the domain of $t$, is $\Sigma_2^{L_{\gamma_0}}(\{T'\})$, in fact its graph is a $\Sigma_1^{L_\alpha}(\{T, \gamma_0\})$ definable subset of $L_{\gamma_0}$ for $\alpha > \gamma_0^+$ the next admissible after $\gamma_0$.

We are thus trying to prove that the starting position $\varnothing$ is good in $L_{\gamma_0}$, and and we have seen that such a quasi-strategy $T^*$ exists satisfying (i) and (ii) in $L_{\gamma_1}$ if and only if such exists in $L_{\gamma_0}$.

Let $B = \cap_{n \in \omega} D_n$ with each $D_n$ recursively open. Define

$$E_n = A \cup \{x \in \lceil T' \rceil \,|\, (\exists p \subseteq x \,(\lceil T'_p \rceil \subseteq D_n \wedge ``p \text{ is not good''}\,)\}.$$

Then, by (2), $E_n$ has a $\Pi_2^{\mathrm{KPI}_0}$ definition (in the parameter $T'$) and again by our elementarity assumptions we have:

$$(4) \;\; x \in L_{\gamma_0} \longrightarrow (\,(x \in E_n)_{L_{\gamma_0}} \longleftrightarrow (x \in E_n)_{L_{\gamma_1}}).$$

Note also that $E_n$ can be considered a $\Sigma_3^0(H)$ set of reals, fixing the good parameter set $H$ as that defined over $L_{\gamma_0}$. Hence, in order to differentiate $E_n$ defined over $L_{\gamma_1}$ (or $L_{\gamma_0}$) and $V$ we set:
$$\tilde{E}_n = A \cup \{x \in \lceil T' \rceil \,|\, (\exists p \subseteq x \,(\lceil T'_p \rceil \subseteq D_n \wedge p \notin H\,)\}.$$

The proof now proceeds by showing

$(+)$        $\exists n \in \omega \; \forall \sigma \in \Sigma_\omega(L_{\gamma_1}) \,(\sigma \text{ is not winning strategy for } I \text{ in } G(\tilde{E}_n;$ $T'))_V$

We first suppose for a moment that this can be shown. We then prove the following:

*Claim* $(+) \Rightarrow \exists T^* \in L_{\gamma_0}(L_{\gamma_0} \vDash ``T^* \text{ witnesses that } \varnothing \text{ is good ''}).$

**Definition 32.** $T''_i =_{\mathrm{df}} \{q \in T' \,|\, \forall p \subseteq q \,(G(E_n, T'_p) \text{ is not a win for } I\,)_{L_{\gamma_i}}\}$ for $i = 0, 1$.

In short, $T''_i$ is $II$'s non-losing quasi-strategy for $G(E_n; T')$, as defined over $L_{\gamma_i}$. ``$(G(E_n, T'_p) \text{ is not a win for } I\,)_{L_{\gamma_i}}$'' means for all $\sigma \in L_{\gamma_i}$ which are strategies for $I$ in $G(E_n, T'_p)$, there is $x \in L_{\gamma_i}$ with $\sigma * x \notin (E_n)_{L_{\gamma_i}}$; this then is a $\Pi_3^{L_{\gamma_i}}$ statement. We may only say that $T''_1 \subseteq T''_0$ since there are potentially more strategies for $I$ in $L_{\gamma_1}$ than in $L_{\gamma_0}$. Neither $T''_i$ is *a priori* a set in $L_{\gamma_i}$, but are definable classes over the respective models. For $i = 0, 1$ we now define $T^*_i \subseteq T''_i$ quasi-strategies for $II$ as follows.

**Definition 33.** $q \in T^*_i$ *if either*:

*(a)* $q \in T_i''$ *and for all positions* $p \subseteq q \lceil T_p' \rceil \nsubseteq D_n$; *or*

*(b) let there be a shortest initial segment* $p \subseteq q$ *with* $p \in T_i''$ *and* $\lceil T_p' \rceil \subseteq D_n$. *In this latter case, by definition of* $T_i''$ $p \in H$. *(The latter must hold, since otherwise we should have that* $(G(E_n; T_p')$ *is a win for* $I)_{L_{\gamma_i}}$ *from position* $p$ *onwards, I making use of some arbitrarily defined but trivial strategy in* $L_{\gamma_0}$, *contradicting that* $p \in T_i''$.) *If the subsequent moves in* $q$ *are consistent with* $t(p) = \hat{T}(p)$, *then we also put* $q$ *into* $T_i^*$. *Otherwise* $q \notin T_i^*$.

(5)   $T_1'' \subseteq T_0''$;   $T_1^* \subseteq T_0^*$;   "$q \in T_1^*$" *is* $\Pi_3^{L_{\gamma_1}}(\{T'\})$.

Note that $T_i^*$ are still quasi-strategies for *II* in $T'$. The following then finishes the Lemma (still under the assumption $( + )$), since as remarked above, "$\varnothing$ is good" is $\Sigma_2^{\mathrm{KPI}_0}(T)$ so the existence of such a witnessing tree $T_0^*$ goes down from $L_{\gamma_1}$ to $L_{\gamma_0}$.

(6)   $(T_0^*$ *witnesses that* $\varnothing$ *is good)*$_{L_{\gamma_1}}$.

Proof: If $x \in \lceil T_0^* \rceil$, then either $x \notin D_n$ or $x \in \lceil \hat{T}(p) \rceil$. In the latter case, as $\hat{T}(p)$ witnesses that $p$ is good, (i) ensures that $x \notin B$. Thus $\lceil T_0^* \rceil \cap B = \varnothing$. We must next show that

*Claim:* $\forall \sigma \in L_{\gamma_1}(\sigma$ *is not a winning strategy for I in* $G(A; T_1^*))$.

Suppose for a contradiction that $\sigma \in L_{\gamma_1}$ were a winning strategy for him in this game.

*Subclaim 1* $\forall y(\sigma * y \in \lceil T_1'' \rceil)$.

Proof: There cannot be a least position $p \in T_1''$ consistent with $\sigma$ so that $\lceil T_p' \rceil \subseteq D_n$: for otherwise for this $p$ we have $T_p^* = \hat{T}(p)$ for the same reason as before: if $p$ were not in $H$ we'd have

$L_{\gamma_1} \vDash$ "$p$ *is not good* $\wedge \lceil T_p' \rceil \subseteq D_n \wedge G(E_n, T_p')$ *is a win for* $I$."

contradicting that $p \in T_1''$.

As there is no such position $p$ like this, we must have that $\forall x$ $\sigma * x \in \lceil T_1'' \rceil$.   QED Subclaim1

*Subclaim 2*  $\forall \sigma_0 \in L_{\gamma_1}(\sigma_0$ *is not a winning strategy for I* in $G(\tilde{E}_n; T_1''))$.

Proof: Suppose for a contradiction $\tau_0 \in L_{\gamma_1}$ were such a winning strategy in $G(\tilde{E}_n \,; T_1'')$.   We convert this to a winning strategy $\tau$ for $I$ in $G(\tilde{E}_n \,, T')$.

We define $\tau$ by initially letting $I$ play using $\tau_0$ until, if ever, $II$ departs from $T_1''$ at some position $p$;  then, as $p \notin T_1''$ we conclude ($I$ has a winning strategy for $G(E_n \,, T_p'))_{L_{\gamma_1}}$ and the latter statement is absolute. $I$ may continue playing using the $L_{\gamma_1}$-least such winning strategy, which we may call $\sigma_p$.

By continuing to play with $\sigma_p$ he wins overall $G(\tilde{E}_n \,, T')$.  The map $\pi$ $p \longmapsto \sigma_p$ is definable over $L_{\gamma_1}$, and hence the overall cumulative strategy $\tau$ we have just implicitly described is also definable over $L_{\gamma_1}$ from $\pi$ and $\tau_0$. It is thus an element of $\Sigma_\omega(L_{\gamma_1})$. If $II$ never departs from $T_1''$ and $I$ uses $\tau_0$ throughout, then this was also a winning run of play in $G(\tilde{E}_n \,, T')$. However the existence of $\tau$ contradicts $(+)$.        QED Subclaim 2

In particular, using Subclaim 2:  as $\sigma \in L_{\gamma_1}$, $\sigma$ itself is not a winning strategy for $I$ in $G(\tilde{E}_n \,; T_1''))$. Hence there is a play $x = \sigma * x_0$ consistent with $\sigma$ satisfying $x \notin \tilde{E}_n$  As $\tilde{E}_n \supseteq A$, we thus have $x \notin A$. However it was originally assumed that $\sigma$ was a winning strategy for $I$ in $G(A; T_1^*)$. By Subclaim 1 $\sigma * x_0 \in [T_1''] \cap [T_1^*]$. This is a contradiction!        QED *Claim*

(7) $(\, G(A; T_0^*) \ is \ not \ a \ win \ for \ I \,)_{L_{\gamma_1}}$.

Proof: By the *Claim*, there is no $\sigma \in L_{\gamma_1}$ a winning strategy for $I$ in $G(A; T_1^*)$.   However $T_1^* \subseteq T_0^*$ and $T_1^*$ does not restrict any of $I$'s moves. Thence (7) holds.   QED (6), (7)

(8) T*here is such a $T^*$ witnessing that $\varnothing$ is good with $T^* \in L_{\gamma_0}$.*

Proof: By (6) and (7) we see that a $T_0^*$ witnessing the requisite $\Sigma_2$ formula can be constructed definably over $L_{\gamma_0}$ from $T$. Thus $T_0^* \in L_{\gamma_1}$. Hence by $\Sigma_2$-reflection of $L_{\gamma_1}$ there is then such a $T^* \in L_{\gamma_0}$.    QED (7)

We thus have to show that (+) above holds. We showed that:

$\varnothing$ is not good $\longrightarrow$
$\forall n \in \omega$ (There is a winning strategy for $I$ in $G(\Delta_4; T')$ definable over $L_{\gamma_1}$).

If we define:

$$\tilde{E}_n^p = A \cup \{x \in \lceil T' \rceil \, | \, (\exists q \subseteq x \,(\, p \subseteq q \wedge \lceil T_q' \rceil \subseteq D_n \wedge q \notin H)\}$$

then the same argument shows that:

(9) $\forall p \in T'($ $p$ is not good $\longrightarrow$ $\forall n \in \omega($There is a winning strategy for $I$ in $G(\bar{E}_n^p; T_p'))$ definable over $L_{\gamma_1}$).

We suppose the lemma false and obtain a contradiction by building a winning strategy $\sigma$ for $I$ for the game $G(A; T')$ which is definable over $L_{\gamma_1}$.

We define the function $s \colon \overline{H} \times \omega \longrightarrow L_{\gamma_1}$ defined by: $s(p, n) = \Sigma_\omega(L_{\gamma_1})$-least winning strategy for $I$ in $G(\bar{E}_n^p; T_p')$. By (9) this function is well defined, and total, on $\overline{H} \times \omega$ and moreover is $\Sigma_\omega^{L_{\gamma_1}}(\{H, T'\})$. Let $\sigma_0 = s(\varnothing, 0)$. Then $\sigma_0$ is a winning strategy for $I$ in $G(\bar{E}_0; T')$. $\sigma$ agrees with $\sigma_0$ until a first, if such occurs, position $p_0$ is reached with $\lceil T_{p_0}' \rceil \subseteq D_0$ but $p_0 \notin H$. If so, then we use the strategy $\sigma_1 = s(p_0, 1)$ for $I$ in $G(\bar{E}_1^{p_0}; T_{p_0}')$. $\sigma$ now agrees with $\sigma_1$ until, if ever, a position $p_2$ is reached with $\lceil T_{p_1}' \rceil \subseteq D_1$ but $p_1$ is not good. The play continues using $\sigma_2 = s(p_2, 2)$. If $q = \cup_{n \in \omega} p_n$ is a non-terminal position, we let $\sigma$ take some arbitrary but canonical choice on positions extending $q$. The strategy $\sigma$ is then definable by a recursion over $J_{\gamma_1+1}$.

(10) If $x \in V$ is any play consistent with $\sigma$ then $x \in A$.

Proof: Suppose firstly we have that for some $n$ $p_n \subseteq x$ is undefined. This implies that there is no initial position $p \subseteq x$ with $(a) p_{n-1} \subseteq p$ (if $n > 0$); $(b)$ $\lceil T_p' \rceil \subseteq D_n$, and $(c)$ $p \notin H$. On the other hand, if all the $p_n$ are defined, then we shall have that $x \in \cap_{n \in \omega} D_n \subseteq B \subseteq A$. Either way we have shown that any play $x \in V$ arising from following the strategy $\sigma$ lies in $A \cap \lceil T' \rceil$. QED(10)

However $T'$ is defined over $L_{\gamma_0}$ to be $II$'s non-losing quasi-strategy in $G(A; T)$ and at (1) it was shown that $T'$ had this property $V$ or equivalently either of the models $L_{\gamma_0}, L_{\gamma_1}$. This contradicts (10)! This finishes the proof of the Lemma.                    QED Lemma

The proof of the theorem now follows Martin [11] pretty much verbatim but again paying attention to definability issues. We repeatedly apply the Lemma with $A = \cup_{n \in \omega} A_n$ and each $A_n \in \Pi_2^0$, acting in turn as an instance of $B$ in the Lemma. This is a $\Sigma_2$-recursion defining a strategy $\tau$ for $II$ over $L_{\gamma_0}$ since all the relevant quasi-strategies given by the Lemma lie in this model. These details now follow.

One applies the lemma with $B = A_0$ obtaining a quasi-strategy for $II$: $T^*(\varnothing)$. By $\Sigma_2$-reflection the $L$-least such lies in $L_{\gamma_0}$, and we shall assume that $T^*(\varnothing)$ refers to it. For any position $p_1 \in T$ with $\mathrm{lh}(p_1) = 1$, let $\tau(p_1)$ be some arbitrary but fixed move in $T'(\varnothing)$, $II$'s non-losing quasi-strategy for the game $G(A, T^*(\varnothing))$. The relation ''$p \in T'(\varnothing)$'' is $\Pi_1^{L_{\gamma_0}}(\{T^*(\varnothing)\})$ and hence ''$y = T'(\varnothing)$'' $\in \Delta_2^{L_{\gamma_0}}(\{T^*(\varnothing)\})$ and thus $T'(\varnothing)$ also lies in $L_{\gamma_0}$. For definiteness we let $\tau(p_1)$ be the numerically least move. For any play, $p_2$ say, of length 2 consistent with the above definition of $\tau$ so far, we apply the lemma again with $B = A_1$ and with $(T^*(\varnothing))_{p_2}$ replacing $T$. This yields a quasi-strategy for $II$, call it $T^*(p_2)$, which is definable in a $\Sigma_2$ way over $L_{\gamma_0}$, in the parameter $(T^*(\varnothing))_{p_2}$. Let $T'(p_2) \in L_{\gamma_0}$ be $II$'s non-losing quasi-strategy for $G(A, T^*(p_2))$, this time with ''$y = T'(p_2)$'' $\in \Delta_2^{L_{\gamma_0}}(\{T^*(p_2)\})$. Again for $p_3 \in T^*(p_2)$ any position of length 3, let $\tau(p_3)$ be some arbitrary but fixed move in $T'(p_2)$. Now we consider appropriate moves $p_4$ of length 4, and reapply the lemma with $B = A_2$ and $(T^*(p_2))_{p_4}$. Continuing in this way we obtain a strategy $\tau$ for $II$, so that $\tau \upharpoonright {}^{2k+1}\omega$, for $k < \omega$, is defined by a recursion that is $\Sigma_2^{L_{\gamma_0}}(\{T\})$. As $L_{\gamma_0} \vDash \Sigma_2$-KP, we have that $\tau \in L_{\gamma_0}$. If $x$ is any play consistent with $\tau$, then for every $n$, by the defining properties of $T^*(p_{2n})$ given by the relevant application of the lemma, $x \in \lceil T^*(x \upharpoonright 2n)\rceil \subseteq \neg A_n$. Hence $x \notin A$, and $\tau$ is a winning strategy for $II$ as required.

<div align="center">QED(<b>Theorem</b>)</div>

**Corollary 34.** *Let $A \in \Sigma_3^0$. Suppose $T$ is a suitable game tree for $A$; then if $(\gamma_0, \gamma_1, \gamma_2)$ are the lexicographic least triple of ordinals with: $\gamma_0 < \gamma_1 < \gamma_2$; with $T \in L_{\gamma_0}$; $L_{\gamma_0} \prec_{\Sigma_2} L_{\gamma_1}$ and $L_{\gamma_0} \prec_{\Sigma_1} L_{\gamma_2} \vDash \mathrm{KPI}_0$, then $L_{\gamma_0} \vDash$ ''$G(A, T)$ is determined''.*